\documentclass{aims}
\usepackage{amsmath}
  \usepackage{paralist}
  \usepackage{graphics} 
  \usepackage{epsfig} 
\usepackage{graphicx}  
\usepackage{epstopdf}
 \usepackage[colorlinks=true]{hyperref}
\hypersetup{urlcolor=blue, citecolor=red}

\def\currentvolume{X}
 \def\currentissue{X}
  \def\currentyear{200X}
   \def\currentmonth{XX}
    \def\ppages{X--XX}
     \def\DOI{10.3934/xx.xx.xx.xx}

\newtheorem{theorem}{Theorem}[section]

\newtheorem{lemma}[theorem]{Lemma}
\newtheorem{corollary}[theorem]{Corollary}

\theoremstyle{definition}

\newcommand{\partder}[2]{\frac{\partial #1}{ \partial #2}}

\newcommand{\D}{\Delta}
\newcommand{\W}{\Omega}

\newcommand{\abs}[1]{\left| #1  \right|}
\newcommand{\infnorm}[1]{\left| \left| #1 \right| \right|_\infty }
\newcommand{\twonorm}[1]{\left| \left| #1 \right| \right|_2}
\newcommand{\heat}{\partial_t - D\Delta}
\newcommand{\heatT}{\partial_t - D_T\Delta}

\newcommand{\R}{\mathbb{R}}

\newcommand*{\defeq}{\mathrel{\vcenter{\baselineskip0.5ex \lineskiplimit0pt
                     \hbox{\scriptsize.}\hbox{\scriptsize.}}}%
     		     =}

\title[Mathematical Analysis of a Spatial Viral Dynamics Model] 
      {Mathematical analysis of an in-host model of viral dynamics with spatial heterogeneity}

\author[Stephen Pankavich and Christian Parkinson]{}


 \email{pankavic@mines.edu}
 \email{christian.a.parkinson@gmail.com}

\thanks{The first author is supported by NSF grant DMS-1211667}

\begin{document}
\maketitle

\centerline{\scshape Stephen Pankavich and Christian Parkinson }
\medskip
{\footnotesize
 \centerline{Colorado School of Mines}
   \centerline{1500 Illinois St.}
   \centerline{Golden, CO 80401, USA}
} 

\medskip


\bigskip

 \centerline{(Communicated by the associate editor name)}

\begin{abstract}
We consider a spatially-heterogeneous generalization of a well-established model for the dynamics of the Human Immunodeficiency Virus-type 1 (HIV) within a susceptible host.  The model consists of a nonlinear system of three coupled reaction-diffusion equations with parameters that may vary spatially. Upon formulating the model, we prove that it preserves the positivity of initial data and construct global-in-time solutions that are both bounded and smooth.  Finally, additional results concerning the local and global asymptotic behavior of these solutions are also provided.

\end{abstract}

\section{Introduction}
Over the past few decades, considerable effort has been devoted to modeling the in-host dynamics of viral infection, and in particular, HIV infection within humans. 
These models have made significant contributions to the understanding of HIV pathogenesis \emph{in vivo} and the mechanisms through which the infection may be mitigated.
The typical response to viral infection within a host is the activation of the immune system, driving the level of virions down. 
If the immune response is sufficiently potent then the disease can be completely eradicated from the body, but often this does not occur. 
Instead, over a time period that can vary from weeks to months, an eventual balance
of viral replication and clearance of the virus by the immune system occurs, leading to a state known as chronic infection. 
These equilibrium outcomes - viral clearance versus the development of a chronic infection - are suggestive
of simple dynamics, but this ignores spatial variations including possible localized regions of infection, which alter the behavior in both transient and long time asymptotics \cite{Haase}. 
 
Although some mathematical modeling of acute HIV infection has been performed \cite{DeB, MKKC, RQC}, many current models focus on the time course of infection during the chronic stage, such as the response to antiretroviral therapy and HIV levels after the viral setpoint has been established \cite{Bonhoeffer_et_al, P2}.
Moreover, because of the intricate complexity and the enormous number of cells involved, mathematical descriptions have generally been limited to nonlinear systems of a few coupled ordinary differential equations describing the average behavior throughout the whole body under the assumption that the environment is well-mixed or spatially-homogeneous \cite{Graw_Perelson,Perelson_Nelson}.
Unfortunately, such an assumption is not valid during the earliest stages of infection or at sites of viral entry, for instance, the vaginal or rectal epithelium during sexual transmission of HIV. 

Indeed, viral propagation is a fundamentally local process. 
Focal bursts of virions have been demonstrated near infected T-cells \emph{in vivo} \cite{Haase, Miller, Reinhart}. 
Correlations between spatial location and viral genetics have been observed and modeled within splenic tissue \cite{Cheynier, Frost}. 
Localization is particularly relevant in the case of HIV infection, for a number of reasons. 
First, the virus is inherently unstable, possibly due to shedding of the glycoprotein gp120, which mediates the virus-cell binding process \cite{McKeating}.
Spatially-heterogeneous outcomes may arise from underlying heterogeneity within tissue architecture, but it is also possible to observe such non-uniformities emerge spontaneously from the infection dynamics. 
Finally, an important aspect of the in-host spread of HIV is that the vast majority of infection occurs only in lymphoid tissues, where target cells are densely packed. 
This particular environment allows the virus to maximize the efficiency of diffusive transfer from infected producers to target cells. 
Thus, the local propagation of HIV within lymphoid tissue is fundamentally different in comparison to the hematogenous spread of the virus to other distant tissues \cite{Grossman, Strain_et_al}.
For all of these reasons, advances in modeling spatial aspects of the development and spread of viral infection within a host are crucial to furthering our understanding of viral pathogenesis and treatment.

Though this paper constitutes the first mathematical study of a model using a system of partial differential equations to study the spatial dependence of infection dynamics, others have used different mathematical mechanisms to model the incorporation of spatial fluctuations and correlations. 
Funk et al \cite{Funk} posited a discrete lattice model to study the spatial dynamics of virions and T-cells, but provided only computational results and did not allow for motion or infection away from a fixed grid.
Strain et al \cite{Strain_et_al} studied the competition between viral lability and diffusion using a spatial cellular automaton model with a fixed T-cell background. 
Finally, Brauner et al \cite{Brauner_et_al} studied a system of coupled ODEs and a single PDE that allowed for the diffusion of virions within a two-dimensional medium under the assumption that T-cells remained motionless.  
However, their study did not focus on the well-posedness or global dynamical behavior of the model, or allow for cellular diffusion.

In order to investigate the impact of spatial dynamics in a simple mathematical model of HIV infection we extend the standard lumped or three-component model of in-host viral dynamics \cite{Bonhoeffer_et_al, Jones_Roemer, Nowak_Bangham, NM, Perelson_Ribeiro} to include spatially random diffusion and a spatially-dependent T-cell supply rate.
Upon describing the model, we establish basic results concerning the well-posedness of smooth solutions and then focus mainly on their longtime asymptotic behavior.
The overall goals are to establish a mathematical framework under which the spatial propagation of infection may be studied and to elucidate the contribution of the parameter space on the dynamical properties of solutions. 

This paper proceeds as follows.  In the next section, we introduce the spatial model of in-host viral dynamics and state our main results concerning existence, uniqueness, regularity, and longtime behavior of solutions to this system of semilinear PDEs.  The proofs of these theorems are contained within Section 3, and finally proofs of associated lemmas are included within Section 4.  Throughout the paper $C$ will denote a positive constant that may change from line to line.  When necessary, we will denote the dependence of this constant on other parameters using subscripts, for instance $C_{t^*}$. 

\section{Spatial Model of Viral Dynamics and Main Results}

\subsection{Derivation of the model}
We begin with the standard, spatially-homogeneous model for in-host virus dynamics, given by 
\begin{equation}
\label{ODE} 
\left. \begin{aligned}
\frac{dT}{dt} &= \lambda - \mu_T T - k TV,\\ 
\frac{dI}{dt}  &= kTV - \mu_I I,\\ 
\frac{dV}{dt} &= N\mu_I I - \mu_V V.
\end{aligned} \right \}
\end{equation} 
Here, the unknown quantities are the populations of (i) uninfected target
cells, denoted by $T(t)$ (in the study of HIV, these are CD4 T-cells); (ii) infected cells, denoted by $I(t)$; and (iii) free virions $V(t)$.
It is assumed within the model that target cells are supplied at a constant rate $\lambda$ and removed either through infection via contact with virions at a rate of $k$ per virion, or through natural cell death with per capita rate $\mu_T$.
Similarly, $\mu_I$ represents the per capita rate at which infected cells are destroyed, either through natural cell death, interaction with the body's immune response, or via lysis due to new virions bursting from the cell membrane, while $\mu_V$ represents the per capita rate at which virions are cleared from the body.
Additionally, $N$ represents the ``burst rate'' of the virus, or the average number of new virions produced over the entire lifespan of an infected cell.

In order to introduce spatial variations within each population, we first let $\Omega \subset \mathbb{R}^n$ be a given bounded domain with $\partial\Omega$ smooth.
As we desire a reformulation that preserves the general biological mechanisms of infection dynamics, we consider the system
\begin{equation}
\tag{3CM}
\label{eq:system} 
\left. \begin{aligned}
\partial_tT - D_T \Delta T &= \lambda (x) - \mu_T T - k TV,\\ 
\partial_tI - D_I \Delta I &= kTV - \mu_I I,\\ 
\partial_tV - D_V \Delta V &= N\mu_I I - \mu_V V
\end{aligned} \right \}
\end{equation} 
for $x \in \Omega$, $t \geq 0$, which describes a spatially-heterogeneous three-component model in which all populations may diffuse throughout a fixed region of the body and the supply rate of new susceptible cells is spatially-dependent, rather than constant.  Hence, $T$, $I$, and $V$ now represent concentrations of healthy cells, infected cells, and virions, respectively.  The diffusive nature of cell transport is, by now, well-known both \emph{in-vivo} and \emph{in-vitro}, and while there may be some minor evidence for this motion to occur via fractional or anomalous diffusion \cite{Harris}, we assume a standard diffusive mechanism via Brownian motion as suggested by \cite{Miller2}.  As spatial effects are present, the domain $\Omega$ might be considered as an entry point of viral infection, such as the vaginal or rectal epithelium during sexual contact.
Additionally, one does not expect the introduction of new T-cells to occur uniformly throughout the domain, and thus, we allow for spatial variations in the function $\lambda(x)$ rather than taking this parameter to be constant as in \eqref{ODE}.

With the spatial domain fixed, the coupled, nonlinear system of PDEs is augmented by the initial conditions
$$ T(0,x) = T_0(x), \quad I(0,x) = I_0(x), \quad V(0,x) = V_0(x).$$
Throughout, we make some standard biological assumptions on the initial concentrations, namely
\begin{equation}
\tag{A1}
\label{A1}
T_0(x), I_0(x), V_0(x) \ \mbox{are continuous, positive, and bounded functions on} \ \bar{\W}
\end{equation}
Additionally, we consider the separate cases of homogeneous Dirichlet or Neumann boundary conditions, so that \eqref{eq:system} is supplemented by either of the two conditions
\begin{equation}
\label{BC1}
\tag{BC1}
w(t,x) = 0 \quad \mathrm{on} \ [0,\infty) \times \partial\Omega
\end{equation}
\begin{equation}
\label{BC2}
\tag{BC2}
\frac{\partial w}{\partial n}(t,x) = 0 \quad \mathrm{on} \ [0,\infty) \times \partial\Omega
\end{equation}
for $w$ replaced by $T, I$, and $V$, respectively.

The biological parameters $k , N, \mu_T, \mu_I, \mu_V$ and diffusion coefficients $D_T, D_I, D_V$ are all positive constants. 
For simplicity, we assume throughout that $\lambda \in C^\infty \left(\Omega \right)$ with $\lambda(x) \geq 0$ for $x \in \Omega$ and $\lambda \not\equiv 0$.

With this framework in place, we make a few remarks concerning the generality of the model.
First, this model can be used to study a wide-range of viral infections, including HIV, Hepatitis B (HBV), Hepatitis C (HCV), Human T-cell Leukemia Virus (HTLV), and Human Cytomegalo Virus (CMV) \cite{Bonhoeffer_et_al, NM}, but we mainly focus on the case of HIV as in-host modeling has significantly furthered the understanding of this particular disease in recent years.
Additionally, though we will deal with the case of $\Omega$ bounded, it is straightforward to use our techniques to establish analogous well-posedness results for the problem posed on the whole space (i.e., $\Omega = \mathbb{R}^n$).
In fact, this problem is easier to approach as one has an exact representation of the associated heat kernel
$$\Phi (t,x) = \frac{1}{(4\pi t)^{n/2}} \exp \left\{ - \frac{\abs{x}^2}{4t}  \right\}, \,\,\, x\in\Omega,\,\, t>0$$ 
and this is the only significant alteration within the associated proofs.
In the case of $\W = \R^n$, one would need to assume $T_0, I_0$, and $V_0$ decay suitably fast as $\vert x \vert \to \infty$, for instance 
$$ T_0(x) \leq C \vert x \vert^{-\alpha}$$
for some $\alpha > n$, $C > 0$, and any $x \in \W$, in lieu of boundary conditions. 
Of course, for this particular application a bounded domain is certainly more natural, and we will typically have in mind the case of Neumann boundary conditions.
Additionally, the cases of $n=1,2,3$ are the most reasonable choices, but our results will actually be independent of dimension.
Hence, we are able to derive sharp conditions on parameters that guarantee the local and global asymptotic stability of solutions regardless of the inherent spatial dimension.
Additionally, though it is assumed throughout that $\lambda(x)$ is smooth, one may relax this assumption and arrive at similar conclusions regarding global-in-time solutions and their dynamical properties.
Finally, another reasonable assumption to include within the model would be to take $D_I = D_T$, as infection of susceptible cells should not influence the rate of diffusion.  However, such a condition will not be needed within the analysis, and hence we omit it.
With the model well formulated, we state the main results of the paper.

\subsection{Existence of Solutions, Positivity, and \emph{a priori} Bounds}
Our first result is quite standard and serves to merely establish the existence of a suitable solution for small time.
\begin{theorem} [Local Existence]
\label{LocalThm}
Assume condition \eqref{A1} holds, then there is $t^* \in [0,\infty]$ such that \eqref{eq:system} with \eqref{BC1}/\eqref{BC2} has a unique, positive solution $(T,I,V) \in C^1((0,t^*); C^2(\Omega)) \times C^1\left ((0,t^*); C^2(\Omega) \right ) \times C^1\left ((0,t^*); C^2(\Omega) \right )$.
\end{theorem}

We will omit the proof as the techniques are well-known in classical literature, and the result can be obtained by a straightforward application of the Contraction Mapping Principle. With a local-in-time solution in hand, we turn our attention to its properties, including positivity (assuming positive initial data), boundedness, and its extension globally in time.

\begin{theorem} [Positivity, Bounds, and Global Existence]
\label{PosBounds}
Assume the initial data satisfy \eqref{A1}. Then, for any $t^* > 0$, there exist $T,I,V$ satisfying \eqref{eq:system} on $(0,t^*) \times \W$, the initial conditions, and boundary condition \eqref{BC1}/\eqref{BC2}. Furthermore, we have $T(t,x),I(t,x),V(t,x) >0$ for all $x \in \W, t \in [0,t^*]$, and there are $C, C_{t^*} > 0$ such that the solutions satisfy 
$$\infnorm{T(t)} \leq C$$
and
$$\infnorm{I(t)} +\infnorm{V(t)} \leq C_{t^*}$$
for all $t \in [0,t^*].$
\end{theorem}

A few comments regarding the global existence of solutions are in order.
Specifically, it should be noted that the idea of ``mass transfer'' among concentrations is a crucial aspect of this model.
Consider, for instance, what would occur if the sign of the $kTV$ term in the first equation were changed, so that $T$ satisfied the evolution equation
$$(\heatT ) T = \lambda (x) - \mu_T T + k TV.$$
Then, even in the simpler case of $\lambda \equiv 0$, it can be shown that this model mimics the behavior of the related nonlinear system
$$(\partial_t u - \D u) =  uv, \qquad (\partial_t v - \D v) =  uv$$
under a suitable transformation of the unknown quantities.
Using tools similar to those established in \cite{EH, Wang}, it can be shown that solutions to this system blow-up in finite time.
Hence, the exact transfer, rather than addition of the nonlinear terms is crucial to the existence theory, as is the continued positivity of solutions.
In addition, Theorem \ref{PosBounds} demonstrates further that the effect of a continuing supply of T-cells, given by $\lambda(x) \geq 0$, does not lead to a blow-up in the system

Now that we have established that solutions exist globally in time, we note that generally, the diffusion operator $( \heat)$ has a smoothing effect on initial data, so we expect some gain in regularity in the concentrations. More specifically, assuming that the initial data is in $L^2(\Omega)$, we expect that solutions not only remain in $L^2(\Omega)$ but actually possess weak derivatives in this space as well.\\

\begin{theorem}[Regularity] 
\label{LOR}
Let $t^*>0$ be given. If $T,I,V$ satisfy \eqref{eq:system} on $(0,t^*) \times \Omega$ and $T_0,I_0,V_0 \in L^2\left(\Omega \right)$ with \eqref{BC1}/\eqref{BC2}, then $\nabla T(t,\cdot), \nabla I(t,\cdot), \nabla V(t,\cdot) \in L^2\left(\Omega \right)$ for all $t\in (0,t^*)$. In fact, $T(t,\cdot), I(t,\cdot), V(t,\cdot) \in H^m\left(\Omega \right)$ for all $t\in (0,t^*)$ and $m \in \mathbb N$. 
\end{theorem}

\vspace{0.1in}

\subsection{Asymptotic Behavior of Solutions}
Finally, with the well-posedness of solutions understood, we turn our attention to their longtime asymptotic behavior. In particular, we study both the local and global stability properties of the uninfected steady state, described below.
Of course, any time-independent solution, i.e. a triple $(T(x), I(x), V(x))$, must satisfy the nonlinear elliptic system 
\begin{equation}
\label{steady} 
\left.
\begin{aligned}
- D_T \Delta T &= \lambda (x) - \mu_T T - kTV\\
- D_I \Delta I &= kTV - \mu_I I\\
- D_V \Delta V &= N\mu_I I - \mu_V V
\end{aligned}
\right \}
\end{equation} 
which may have many solutions for differing parameter values.
In particular, we focus on the unique uninfected state $E_c := (T_\infty(x), 0, 0)$, which represents the ultimate clearance of the virus, where $T_\infty$ satisfies the linear boundary value problem
\begin{equation}
\label{Tinf} 
- D_T \Delta T_\infty = \lambda (x) - \mu_T T_\infty
\end{equation} 
for $x \in \Omega$ with boundary conditions given by \eqref{BC1}/\eqref{BC2}.
We first summarize some useful properties of the steady-state T-cell count.
\begin{theorem}
\label{Tinfprop}
The function $T_\infty$ defined by \eqref{Tinf} with \eqref{BC1}/\eqref{BC2} is $C^\infty$ and satisfies
$$0 \leq T_\infty(x) \leq \frac{\Vert \lambda \Vert_\infty}{\mu_T}$$
for every $x \in \Omega$ with $T_\infty \not\equiv 0$.
\end{theorem}

In order to determine the stability of the aforementioned equilibrium state, we are first led to study the system \eqref{steady} linearized about $E_c$, namely 
\begin{equation}
\label{lin} 
\left.
\begin{aligned}
- D_T \Delta T &= \lambda (x) - \mu_T T - kT_\infty(x)V\\
- D_I \Delta I &= kT_\infty(x)V - \mu_I I\\
- D_V \Delta V &= N\mu_I I - \mu_V V
\end{aligned}
\right \}
\end{equation} 
with boundary conditions \eqref{BC1}/\eqref{BC2}.
Here, the first equation decouples from the last two and we need only study the PDEs describing $I(x)$ and $V(x)$. Accordingly, we define the second-order, linear self-adjoint operator
\begin{equation}
\label{L}
\mathcal{L} = \nabla \cdot (D \nabla) + M(x)
\end{equation}
on the Hilbert space
$$\mathcal{H} = \left \{\phi \in H^2(\Omega) \times H^2(\Omega) : \phi_1 \ \mathrm{and} \  \phi_2 \ \mathrm{satisfy} \ \eqref{BC1}/ \eqref{BC2} \right \}$$
where
\begin{equation}
\label{DM}
D = \mathrm{diag}(D_I, D_V) \qquad \mathrm{and} \qquad M(x) = \left[ \begin{array}{cc}
-\mu_I & kT_\infty(x) \\
N\mu_I & -\mu_V  \end{array} \right].
\end{equation}
We note that $\mathcal{L}$ is a negative operator if $-M(x)$ is a positive semi-definite matrix.  Because $\mu_I > 0$ this latter condition is implied by $\det(-M(x)) > 0$ for all $x \in \Omega$. After some algebra, this condition is exactly $\Vert R_0 \Vert_\infty < 1$ where
\begin{equation}
\label{R0x}
R_0(x) := \frac{NkT_\infty(x)}{\mu_V}.
\end{equation}

The stability properties of $E_c$ then depend crucially on the greatest real part of the spectrum of $\mathcal{L}$, namely
\begin{equation}
\label{eta0}
\eta_0 := \sup \{ \mathrm{Re}(\eta) : \eta \in \sigma(A) \}.
\end{equation}
Due to the compactness of the inverse of $\mathcal{L}$, it is well-known that $\sigma(\mathcal{L})$ consists only of eigenvalues and, when ordered, these form a decreasing sequence of real numbers diverging to $-\infty$.
Additionally, $\eta_0$ can be uniquely identified using the associated Rayleigh quotient
$$\eta_0 =  \sup_{\substack{\psi \in \mathcal{H}\\ \psi \neq 0}} \left \{ \frac{\int_\Omega \psi(x) \cdot \mathcal{L} \psi(x) \ dx}{\int_\Omega \vert \psi(x) \vert^2 \ dx} \right \}.$$
Hence, the sign of the parameter $\eta_0$ will determine whether all eigenvalues are negative, ensuring that solutions which begin sufficiently close tend to the equilibrium state, or whether a positive eigenvalue exists, thereby creating instabilities within the system.
Prior to stating the stability and instability results for the viral clearance state, we first show that under particular conditions it is the unique biologically-feasible solution of the steady state system \eqref{steady}.
\begin{theorem}
\label{unique}
If $\eta_0 < 0$ then the triple $(T_\infty(x), 0, 0)$ is the unique nonnegative solution of \eqref{steady} with \eqref{BC1}/\eqref{BC2}.
\end{theorem}

Next, we characterize the local stability properties of the clearance state.  
In the spatially-homogeneous case, as in many other biological systems, it is known that a single parameter dictates the dynamical structure of solutions.  With regards to \eqref{ODE}, the quantity known as the basic reproductive ratio, defined by $R_0 := \frac{\lambda k N}{\mu_T \mu_V}$, describes the entire long time behavior of the infection.  In this case the spatially-independent T-cell count is $T_\infty = \frac{\lambda}{\mu_T}$, which means that $R_0(x)$ as defined by \eqref{R0x} is exactly the basic reproductive ratio.  Additionally, the associated principal eigenvalue is $\eta_0 = R_0 - 1$. Though $R_0$ for the spatially-heterogeneous system is now a function rather than a single value, we may still characterize the dynamics in terms of the single parameter $\eta_0$.
Hence, the following local stability theorem establishes a parameter regime that generalizes the spatially-homogeneous model, in which the longtime asymptotics are still shown to be determined only by the sign of $\eta_0$.

\begin{theorem}
\label{local}
Let $T,I,V$ satisfy \eqref{eq:system} with \eqref{BC1}/\eqref{BC2}.
If $\eta_0 < 0$ then the viral clearance state $E_c$ is locally asymptotically stable.  If $\eta_0 > 0$, then it is unstable.
\end{theorem}

Finally, we study the global dynamics of this steady state. 
Our first global result shows that the T-cell population tends to $T_\infty(x)$ exponentially fast for any initial data and within any parameter regime.
\begin{theorem}
\label{Tinfty}
Let $T,I,V$ satisfy \eqref{eq:system} with \eqref{BC1}/\eqref{BC2}. Then, for every $ t \geq 0$
$$ \Vert T(t) - T_\infty \Vert_\infty \leq \Vert T_0 - T_\infty \Vert_\infty e^{-\mu_T t}.$$
\end{theorem}
\noindent In view of Theorems \ref{unique} and \ref{local}, we see that for $\eta_0 < 0$, $E_c$ is locally stable and the only nonnegative equilibrium.  Hence, one might expect that it is also globally stable in this case.  While we do not have a proof of this result assuming only the condition $\eta_0 < 0$, we can do so under a slightly more restrictive assumption.
\begin{theorem}
\label{Asymp}
Let $T,I,V$ satisfy \eqref{eq:system} with \eqref{BC1}/\eqref{BC2}.
If $\Vert R_0 \Vert_\infty < 1$ then there are $C_0, a > 0$ such that
$$\Vert I(t) \Vert_\infty + \Vert V(t) \Vert_\infty \leq C_0e^{-at}$$
for every $t \geq 0$.
\end{theorem}
\noindent Combining this result with Theorem \ref{Tinfty} yields a sufficient condition under which the viral clearance state is a globally asymptotically stable equilibrium point of the system with exponential decay.
Though it is beneficial to characterize the dynamics of solutions, it may be difficult to explicitly determine $T_\infty$ for given parameter values and domain $\Omega$.  Hence, our final result provides an additional sufficient condition that can be readily computed to describe the global asymptotics of solutions. 
 
\begin{corollary}
\label{global}
Let $T,I,V$ satisfy \eqref{eq:system} with \eqref{BC1}/\eqref{BC2}.  If 
$\frac{Nk \Vert \lambda \Vert_\infty}{\mu_T \mu_V} < 1,$
then the conclusion of Theorem \ref{Asymp} remains valid.
\end{corollary}

We note that the longtime behavior results are independent of both the diffusion coefficients $D_T, D_I, D_V > 0$ and the dimension $n = 1, 2, 3$. 
In the next section we prove the main results of the paper.

\section{Proofs of Main Results}
As previously mentioned, the proof of the first theorem is a fairly standard application of the Contraction Mapping Principle, and hence we omit it.
Prior to proving Theorem \ref{PosBounds}, we establish a few minor lemmas regarding the scalar diffusion equation that arise from the Maximum Principle.

\begin{lemma}
\label{Pos}
Let $D > 0$ and $u_0(x) > 0$ for all $x\in\W$ be given and assume that $u$ satisfies
\begin{equation} 
 \label{eq:scal_heat}
\left. \begin{array}{lr} &(\heat ) u = g(t,x), \,\,\,\,\,\,\,\,\, x \in \W, t>0 \\ &u(0,x) = u_0(x). \hspace{1.4cm} x \in \W \end{array} \right \}
\end{equation}
If $g(t,x) \ge 0$ for all $x \in \W, t \geq 0$, then $u(t,x) \geq 0$ for all $x \in \W, t \geq 0$.\\  
Additionally, if $g(t,x) > 0$ for all $x \in \W, t \geq 0$, then $u(t,x) > 0$ for all $x \in \W, t \geq 0$.
That is, $(\heat)$ preserves positivity on $[0,\infty) \times \W$.
\end{lemma}

The following corollaries, which will be useful throughout the remainder of this section, are then immediate consequences of this result.

\begin{corollary} 
\label{heatbound1}
If $u(t,x)$ satisfies \eqref{eq:scal_heat} with $D > 0$, then $$\infnorm{u(t)} \le \infnorm{u_0} + \int^t_0 \infnorm{g(\tau)} d\tau, \,\,\, \text{ for all } t \in [0,\infty).$$ 
\end{corollary}

\begin{corollary} 
\label{heatbound2}
Assume that $u(t,x)$ satisfies the differential inequality $$\begin{aligned} &(\heat ) u \le g(t,x), \,\,\,\,\,\,\,\,\, x \in \W, t > 0 \\ &u(x,0) = u_0(x), \hspace{1.4cm} x \in \Omega. \end{aligned}$$ Then $u(t,x)$ satisfies the same inequality as in Corollary \ref{heatbound1}. That is, $$\infnorm{u(t)} \le \infnorm{u_0} + \int^t_0 \infnorm{g(\tau)} d\tau, \,\,\, \text{ for all } t \in [0,\infty).$$
\end{corollary}

With these results in place, we may now prove the global existence and positivity theorem.

\begin{proof}[Proof (Theorem \ref{PosBounds}).]

Since the initial data is positive and solutions are smooth (see Theorem \ref{LOR}), the functions $T(t,x), I(t,x), V(t,x)$ must remain positive for all $x \in \W$ and up to some time within the interval of existence $[0,t^*).$ 
Define the maximal time of positivity within this interval by
$$ Q= \sup \left\{ t \in [0,t^*) \,\, \biggr \vert \,\,  T(s,x), I(s,x), V(s,x) >0  \text{ for all } x \in \W,  s\in [0,t] \right\}.$$ 
Estimating on this interval we find
$$(\partial_t - D_T\D) T + \mu_TT \le \lambda(x). $$ 
Multiplying by an integrating factor yields 
$$(\partial_t - D_T\D) [e^{\mu_T t} T] \le \lambda(x) e^{\mu_T t}. $$ 
By Corollary \ref{heatbound2}, we find $$\infnorm{T(t)e^{\mu_Tt}} \le \infnorm{T_0}+ \int^t_0 \infnorm{\lambda} e^{\mu_Ts} ds$$ 
and upon integrating
$$\infnorm{T(t)} \le \infnorm{T_0}e^{-\mu_T t} + \frac{\infnorm{\lambda}}{\mu_T} \left( 1 - e^{-\mu_Tt} \right).$$ 
In particular, we have a uniform bound for $\infnorm{T(t)}$ on $[0,Q]$, denoted by 
$$T_M: = \infnorm{T_0} + \frac{\infnorm{\lambda}}{\mu_T}.$$

Next, we use the positivity of solutions and this uniform bound within the equations for $I$ and $V$.
On the interval $[0,Q]$, we have $\left (\partial_t - D_I\D \right) I \le kT(t,x)V(t,x)$. 
Using Corollary \ref{heatbound2}, this implies 
\begin{equation}
\begin{aligned}[b] 
\label{eq:Isup_bound}
\infnorm{I(t)} &\le\infnorm{I_0} + k \int^t_0 \infnorm{T(\tau)V(\tau)} d\tau \\ 
& \le \infnorm{I_0} + k T_M \int^t_0 \infnorm{V(\tau)} d\tau \\ 
& \le C \left( 1+  \int^t_0 \infnorm{V(\tau)} d\tau\right). 
\end{aligned}
\end{equation}
Using the positivity of $I$ within the equation for $V$, we have
$$(\partial_t - D_V\D)V  \le N\mu_I I(t,x)$$ for $t\in [0,Q].$ 
Invoking Corollary \ref{heatbound2} yields
\begin{equation} \begin{aligned}[b] \label{eq:Vsup_bound} \infnorm{V(t)} &\le \infnorm{V_0} + N\mu_I \int^t_0 \infnorm{I(\tau)} d\tau \\ & \le C\left( 1 + \int^t_0 \infnorm{I(\tau)} d\tau \right). \end{aligned}\end{equation} 
Finally, define $$\phi(t) = \infnorm{I(t)} + \infnorm{V(t)}, \,\,\,\, t \in [0,Q].$$ 
Adding \eqref{eq:Isup_bound} and \eqref{eq:Vsup_bound}, we see 
$$\phi(t) \le C \left( 1 + \int^t_0 \phi(\tau) d\tau\right), \,\,\,\, t \in [0,Q].$$ 
By Gronwall's Inequality, we can conclude that $$\phi(t) \le  C e^{Ct}, \,\,\,\, t\in [0,Q].$$ 
Thus, both $\infnorm{I(t)}$ and $\infnorm{V(t)}$ remain finite on $[0,Q]$.  
Therefore, solutions can be continued indefinitely, as long as they remain positive.

Now, these bounds will imply the continuing positivity of solutions to time $t = Q$.
Rearranging terms within the equation for $I$ yields 
$$(\partial_t - D_I\D) I + \mu_I I= kT(t,x)V(t,x) >0, \,\,\,\, t\in [0,Q].$$ 
Multiplying this equation by an integrating factor, we arrive at 
$$(\partial_t - D_I\D ) \left [e^{\mu_I t} I \right] = k e^{\mu_I t} T(t,x)V(t,x).$$  
Thus, we have a diffusion equation for $e^{\mu_I t}I(t,x)$ with a nonnegative forcing term and positive initial conditions. By Lemma \ref{Pos}, we can conclude that $e^{\mu_I t} I(t,x) >0$ for all $t\in [0,Q]$, and hence $I(t,x)>0$ for $t \in [0,Q]$.
Positivity of $V$ follows in the same way since $V$ satisfies 
$$(\partial_t - D_V\D) V + \mu_V V = N\mu_I I(t,x).$$ 
and this becomes
$$(\partial_t - D_V\D)\left [ e^{\mu_V t} V\right ] = N\mu_I e^{\mu_V t} I(t,x).$$ 
Thus, by Lemma \ref{Pos}, $e^{\mu_V t} V(x,t)>0$ for $t \in [0,Q]$ and $V(t,x) >0$ for $t\in [0,Q]$.

Finally, for the $T$ equation, we use the local-in-time bound on $\infnorm{V(t)}$ to find 
\begin{align*} 
(\partial_t - D_T\D) T &\ge \lambda - \mu_TT - k T(t,x) \infnorm{V(t)} \\ 
& \ge \lambda - \mu_TT - Ce^{Ct} T
\end{align*} 
Rearranging yields $$(\partial_t - D_T\D) T + \left[ \mu_T + Ce^{Ct}\right] T \ge \lambda$$ 
and using the integrating factor $\zeta(t) = \exp \left\{ \int_0^t (\mu_T + Ce^{C\tau}) \, d\tau\right\}$, we find
$$ (\partial_t - D_T\D)  \left\{ \zeta(t) T \right\} \ge \lambda \zeta(t)$$ 
and $\zeta(0)T(0,x) = T_0(x)$. 
Thus, by Lemma \ref{Pos}, we have $\zeta(t)T(x,t) >0$, and hence $T(t,x)>0$ for $t\in [0,Q]$. 
Therefore, solutions remain strictly positive throughout the interval $[0,Q]$ and this implies $Q = t^*$, the maximal time of existence.  Of course, since the solution remains bounded on this interval as well, we find $t^* = \infty$ and the proof is complete.
\end{proof}

Next, we turn to the proofs of the regularity results.  For brevity, we will prove the first statement only. However, the same technique can be applied to derive estimates inductively and prove the latter result without the introduction of new ideas.  For additional background on the specifics of obtaining the second statement, we direct the reader to \cite{Pankavich_Michalowski1, Pankavich_Michalowski2} as the higher-order regularity stated here can be deduced by straightforwardly adapting the arguments in those previous works.

\begin{proof}[Proof (Theorem \ref{LOR})]
Beginning with the equation for $T$, we multiply by $T$ and integrate over $\Omega$ to arrive at 
$$ \frac{1}{2}\frac{d}{dt}\twonorm{T(t)}^2 - D_T \int_{\W} T \D T dx = \int_{\W} \lambda T dx - \mu_T \twonorm{T(t)}^2 - k\int_{\W} T^2V dx. $$ 
Integrating by parts on the left, enforcing boundary conditions, using the Cauchy-Schwarz inequality, and replacing $V$ by it's supremum, we see 
\begin{align*}\frac{1}{2} \frac{d}{dt}\twonorm{T(t)}^2 & + D_T \twonorm{\nabla T(t)}^2\\ &\le \twonorm{\lambda} \twonorm{T(t)} - \mu_T \twonorm{T(t)}^2 + k \infnorm{V(t)}\twonorm{T(t)}^2 \\ &\le \tfrac{1}{2}\left(\twonorm{\lambda}^2 + \twonorm{T(t)}^2 \right) -\mu_T \twonorm{T(t)}^2 + k \infnorm{V(t)}\twonorm{T(t)}^2. \end{align*}
Finally we arrive at 
\begin{equation} \label{eq:dTdtbound} \frac{d}{dt}\twonorm{T(t)}^2 \le C_{t^*}\left( 1 + \twonorm{T(t)}^2 \right) - 2D_T \twonorm{\nabla T(t)}^2.
\end{equation} 
To generate a derivative estimate, we proceed similarly and take the gradient of the equation dotted with $\nabla T$ so that 
$$\frac{1}{2} \partial_t \left|\nabla T\right|^2 - D_T \nabla T \cdot \nabla \D T = \nabla T \cdot \nabla \lambda - \mu_T \left| \nabla T \right|^2 - k \nabla T \cdot \left( V\nabla T + T \nabla V \right).$$ 
Integrating over spatial variables and using integration by parts yields
\begin{align*} & \frac{1}{2} \frac{d}{dt} \twonorm{\nabla T(t)}^2 + D_T \twonorm{\D T(t)}^2\\ 
& \quad= \int_{\W} \nabla T \cdot \nabla \lambda dx - \mu_T \twonorm{\nabla T(t)}^2 - k \int_{\W} \left( V \left| \nabla T \right|^2 + T \nabla T \cdot \nabla V \right)  dx \\
& \quad \le \tfrac{1}{2} \left(\twonorm{\nabla T(t)}^2 + \twonorm{\nabla \lambda}^2 \right) -\mu_T \twonorm{\nabla T(t)}^2 + k \left( \infnorm{V(t)}\twonorm{\nabla T(t)}^2 +\right. \\ &\hspace{1cm}+\left. \tfrac{1}{2}\infnorm{T(t)} \left( \twonorm{\nabla T(t)}^2 + \twonorm{\nabla V(t)}^2\right) \right), 
\end{align*} 
whence 
\begin{equation} \label{eq:gradTbound} \frac{d}{dt} \twonorm{\nabla T(t)}^2 \le C_{t^*} \left( 1+ \twonorm{\nabla T(t)}^2 + \twonorm{V(t)}^2 \right) - 2 D_T \twonorm{\D T(t)}^2. 
\end{equation} 
We deal with $I$ analogously. Multiplying this evolution equation by $I$ and integrating gives 
\begin{align*} \frac{1}{2} \frac{d}{dt} \twonorm{I(t)}^2 + D_I \twonorm{\nabla I(t)}^2 &= k \int_{\W} TIV dx - \mu_I \twonorm{I(t)}^2 \\ &\le \infnorm{V(t)} \int_{\W} TI dx - \mu_I \twonorm{I(t)}^2 \\ & \le \tfrac{1}{2}\infnorm{V(t)} \left( \twonorm{T(t)}^2 + \twonorm{I(t)}^2 \right) -\mu_I \twonorm{I(t)}^2. 
\end{align*}
From this we find
\begin{equation}\label{eq:dIdtbound} \frac{d}{dt} \twonorm{I(t)}^2 \le C_{t^*} \left( \twonorm{T(t)}^2 + \twonorm{I(t)}^2 \right) - 2 D_I \twonorm{\nabla I(t)}^2. 
\end{equation} 
Next, taking the gradient of this equation and then dotting with $\nabla I$, we arrive at
\begin{align*} \frac{1}{2} \partial_t \abs{\nabla I}^2 - D_I \nabla I \cdot \nabla \D I = k \nabla I \cdot \left( V\nabla T + T \nabla V\right) - \mu_I \abs{\nabla I}^2.
\end{align*} 
Integrating over $\Omega$ yields 
\begin{align*}& \frac{1}{2}\frac{d}{dt} \twonorm{\nabla I(t)}^2 + D_I \twonorm{\D I(t)}^2 + \mu_I \twonorm{\nabla I(t)}^2\\  & \quad = k \left( \int_{\W} V \nabla I \cdot \nabla T dx + \int_{\W} T \nabla I \cdot \nabla V dx \right)\\ & \quad \le k\left( \infnorm{V(t)} \int_{\W} \nabla I \cdot \nabla T dx + \infnorm{T(t)} \int_{\W} \nabla I \cdot \nabla V dx \right)\\ & \quad\le C_{t^*} \bigg( \tfrac{1}{2} \left( \twonorm{\nabla I(t)}^2 + \twonorm{\nabla T(t)}^2 \right) +\tfrac{1}{2} \left( \twonorm{\nabla I(t)}^2 + \twonorm{\nabla V(t)}^2 \right) \bigg).
\end{align*} 
and this produces the inequality 
\begin{equation} \label{eq:gradIbound}\frac{d}{dt} \twonorm{\nabla I(t)}^2 \le C_{t^*} \left( \twonorm{\nabla T(t)}^2 + \twonorm{\nabla I(t)}^2  + \twonorm{\nabla V(t)}^2 \right) - 2D_I\twonorm{\D I(t)}^2. 
\end{equation} 
Within the $V$ equation we multiply by $V$ and integrate to find 
\begin{align*} \frac{1}{2} \frac{d}{dt} \twonorm{V(t)}^2 + D_V \twonorm{\nabla V(t)}^2 &= N\mu_I \int_{\W} IV dx - \mu_V \twonorm{V(t)}^2 \\ &\le \frac{N\mu_I}{2}\left( \twonorm{I(t)}^2 + \twonorm{V(t)}^2\right) - \mu_V \twonorm{V(t)}^2 \\ &\le C_{t^*} \left( \twonorm{I(t)}^2 + \twonorm{V(t)}^2 \right). 
\end{align*} 
This yields 
\begin{equation} \label{eq:dVdtbound} \frac{d}{dt} \twonorm{V(t)}^2 \le C_{t^*} \left( \twonorm{I(t)}^2 + \twonorm{V(t)}^2 \right) - 2D_V \twonorm{\nabla V(t)}^2. 
\end{equation} 
The analogous derivative estimate is then
\begin{align*} \frac{1}{2} \frac{d}{dt} \twonorm{\nabla V(t)}^2 + D_V \twonorm{\D V(t)}^2 &= N\mu_I \int_{\W} \nabla I \cdot \nabla V dx - \mu_V \twonorm{\nabla V(t)}^2 \\ &\le \frac{N\mu_I}{2}\left( \twonorm{\nabla I(t)}^2 + \twonorm{\nabla V(t)}^2\right) - \mu_V \twonorm{\nabla V(t)}^2 \\ & \le C_{t^*} \left( \twonorm{\nabla I(t)}^2 + \twonorm{\nabla V(t)}^2 \right), \end{align*} 
from which it follows that 
\begin{equation} \label{eq:gradVbound}\frac{d}{dt} \twonorm{\nabla V(t)}^2 \le C_{t^*} \left( \twonorm{\nabla I(t)}^2 + \twonorm{\nabla V(t)}^2 \right) - 2D_V \twonorm{\D V(t)}^2. \end{equation}

Finally, let $D_{\min} = \min\{D_T,D_I,D_V \}$ and for $t\in [0,t^*]$, define \begin{align*} M(t) &= \Big( \twonorm{T(t)}^2+ \twonorm{I(t)}^2+ \twonorm{V(t)}^2\Big)\\
&  \qquad + D_{\min} t \Big(\twonorm{\nabla T(t)}^2+\twonorm{\nabla I(t)}^2+\twonorm{\nabla V(t)}^2 \Big) \\ &=: \phi_0(t) + t\phi_1(t). \end{align*} 

By adding equations \eqref{eq:dTdtbound},\eqref{eq:dIdtbound} and \eqref{eq:dVdtbound}, we see $$\phi_0'(t) \le C_{t^*}\left( 1 + \phi_0(t) \right) - 2D_{\min} \phi_1(t)$$ and by adding equations \eqref{eq:gradTbound}, \eqref{eq:gradIbound} and \eqref{eq:gradVbound}, we see $$\phi'_1(t) \le C_{t^*}(1+\phi_1(t)) - 2D_{\min} \Big(\twonorm{\D T(t)}^2+ \twonorm{\D I(t)}^2+ \twonorm{\D V(t)}^2 \Big).$$ 
Then, letting $\phi_2(t) = \twonorm{\D T(t)}^2+ \twonorm{\D I(t)}^2+ \twonorm{\D V(t)}^2,$ we arrive at 
\begin{align*} M'(t) &= \phi'_0(t) + D_{\min}\phi_1(t) + D_{\min}t \phi'_1(t) \\ &\le C_{t^*}\left( 1 + \phi_0(t) \right) - 2D_{\min}\phi_1(t) + D_{\min}\phi_1(t) +D_{\min}t\Big(C_{t^*}(1+\phi_1(t)) - 2\phi_2(t) \Big) \\ &\le C_{t^*}\big(1+\phi_0(t) + D_{\min}t\phi_1(t) \big) - D_{\min}\phi_1(t) - 2D_{\min}t\phi_2(t). \end{align*} 
Noting that $\phi_1$ and $\phi_2$ are nonnegative, we find 
$$M'(t) \le C_{t^*} \big( 1+M(t) \big).$$ 
From this, an application of Gronwall's inequality produces the bound 
$$M(t) \le C_{t^*}\big(1 + M(0) e^t \big) \le C_{t^*} (1+M(0)).$$ 
By assumption $$M(0) = \twonorm{T_0}^2 + \twonorm{I_0}^2 + \twonorm{V_0}^2$$ is finite and thus $M(t)$ remains finite on the interval.
Since $\phi_0$ is nonnegative, this implies the bound 
$$\phi_1(t) \le \frac{C_{t^*}}{D_{\min}t}(1+M(0))$$
for $t\in (0,t^*].$
As each of the quantities $\twonorm{\nabla T(t)}^2, \twonorm{\nabla I(t)}^2, \twonorm{\nabla V(t)}^2$ is bounded by $\phi_1(t)$, we find $$\nabla T(t, \cdot), \,\, \nabla I(t, \cdot), \,\, \nabla V(t,\cdot) \in L^2(\Omega)$$ for any $t \in (0,t^*]$.
To prove the latter statement of the theorem, higher-order estimates are obtained and combined using the same techniques.  The result then follows by induction (see \cite{Pankavich_Michalowski2} for further details).
\end{proof}

With the regularity proof complete, we finally study the longtime dynamical behavior of \eqref{eq:system}, and in particular, the viral clearance steady state $E_c$.  We first prove some useful properties of $T_\infty(x)$.

\begin{proof}[Proof (Theorem \ref{Tinfprop})]
To prove the nonnegativity of $T_\infty$, we merely note that the Green's function associated to the operator $(-D_T\Delta + \mu_T)$ is positive \cite{Kreith}, and hence the conclusion follows due to the nonnegativity of $\lambda(x)$ on $\Omega$.  Additionally, $T_\infty$ inherits the regularity of $\lambda$ and $T_\infty \not\equiv 0$ because $G,\lambda \not \equiv 0$.  

To prove the upper bound, we let $G(x)$ satisfy
$$(-D_T\Delta + \mu_T) G = \delta(x)$$
for $x \in \Omega$.  Then, integrating over $\Omega$, we find
$$\int_\Omega G(x) \ dx  = \frac{1}{\mu_T} + \frac{D_T}{\mu_T} \int_{\partial \Omega} \frac{\partial G}{\partial n} \ dS.$$
If $G$ satisfies \eqref{BC2}, we merely enforce the boundary conditions to conclude
$$\int_\Omega G(x) \ dx  = \frac{1}{\mu_T}.$$
If $G$ satisfies \eqref{BC1}, we note that $G = 0$ on $\partial \Omega$ and $G(x) > 0$ for $x \in \Omega$, therefore $\frac{\partial G}{\partial n} \leq 0$ on $\partial \Omega$ and thus
$$\int_\Omega G(x) \ dx  \leq \frac{1}{\mu_T}.$$
Again, since $G$ is positive, this further shows $\Vert G \Vert_1 \leq \frac{1}{\mu_T}$ in either case.
Finally, since $T_\infty$ satisfies the inhomogeneous equation it must be given by $T_\infty = G \ast \lambda$ and a simple convolution estimate yields
$$\Vert T_\infty \Vert_\infty = \Vert G \ast \lambda \Vert_\infty \leq \Vert \lambda \Vert_\infty \Vert G \Vert_1 = \frac{\Vert \lambda \Vert_\infty}{\mu_T}.$$
\end{proof}

Next, we prove the uniqueness of this steady state when $\eta_0 < 0$.

\begin{proof}[Proof (Theorem \ref{unique})]
Let $(\tilde{T}(x), \tilde{I}(x), \tilde{V}(x)) \neq (T_\infty(x), 0, 0)$ be another nonnegative solution of \eqref{steady} satisfying \eqref{BC1}/\eqref{BC2}.
In particular, we must have $\tilde{I} \not\equiv 0$ and $\tilde{V} \not\equiv 0$ since either condition implies the other and $\tilde{T} \equiv T_\infty$.
The second and third equations in the system, namely
$$\begin{aligned}
- D_I \Delta \tilde{I} &= k\tilde{T} \tilde{V} - \mu_I \tilde{I}\\
- D_V \Delta \tilde{V} &= N\mu_I \tilde{I} - \mu_V \tilde{V}
\end{aligned}$$
can be rewritten as
$$\begin{aligned}
D_I \Delta \tilde{I} + kT_\infty \tilde{V} - \mu_I \tilde{I} &= k(T_\infty - \tilde{T}) \tilde{V}\\
D_V \Delta \tilde{V} + N\mu_I \tilde{I} - \mu_V \tilde{V} &= 0
\end{aligned}$$
or
\begin{equation}
\label{operator}
\mathcal{L} \left [\begin{array}{c} \tilde{I}\\ \tilde{V} \end{array} \right ] = \left [\begin{array}{c} k(T_\infty - \tilde{T}) \tilde{V}\\ 0 \end{array} \right ].
\end{equation}

Next, let $u(x) = \tilde{T}(x) - T_\infty(x)$ where $T_\infty$ is defined as the unique solution of the linear equation \eqref{Tinf} with \eqref{BC1}/\eqref{BC2}.
Then, $u$ satisfies
$$\left [ -D_T \Delta + \mu_T \right ] u = - k\tilde{T}\tilde{V}.$$
Since $\tilde{T}(x), \tilde{V}(x) \geq 0$ for all $x \in \Omega$, it follows that $u(x) \leq 0$ for all $x \in \Omega$.
Therefore, equation \eqref{operator} can be rewritten as
$$ 
\mathcal{L} \left [\begin{array}{c} \tilde{I}\\ \tilde{V} \end{array} \right ] = \left [\begin{array}{c} -ku\tilde{V}\\ 0 \end{array} \right ]$$
Taking the dot product of this equation with $ \left [\begin{array}{c} \tilde{I}\\ \tilde{V} \end{array} \right ]$ and integrating over $\Omega$, we find
$$\int_\Omega  \left [\begin{array}{c} \tilde{I}(x)\\ \tilde{V}(x) \end{array} \right ] \cdot \mathcal{L}  \left [\begin{array}{c} \tilde{I}(x)\\ \tilde{V}(x) \end{array} \right ] \ dx = k\int_\Omega -u(x) \tilde{I}(x) \tilde{V}(x) \ dx \geq 0. 
$$
However, since $\eta_0 < 0$ and $\tilde{I}, \tilde{V} \not\equiv 0$, we find
$$\int_\Omega  \left [\begin{array}{c} \tilde{I}(x)\\ \tilde{V}(x) \end{array} \right ] \cdot \mathcal{L}  \left [\begin{array}{c} \tilde{I}(x)\\ \tilde{V}(x) \end{array} \right ] \ dx
\leq \eta_0 \left ( \Vert \tilde{I} \Vert_2 + \Vert \tilde{V} \Vert_2 \right ) < 0.$$
This inequality contradicts the previous one, and hence the original assumption of an additional solution cannot hold.
\end{proof}

Next, we prove the local stability and instability results.
Since a generalization of the celebrated next generation method adapted to spatially-heterogeneous problems was recently developed in \cite{WangZhao}, we will rely on many of the results therein for the proof.
\begin{proof}[Proof (Theorem \ref{local})]
To prove the first conclusion, we assume $\eta_0 < 0$ and put the system \eqref{lin} into the form of \cite[Equation (3.9)]{WangZhao}.
The uninfected compartment - $T(t,x)$ - and infected compartments - $I(t,x)$ and $V(t,x)$ - can be separated into the scalar function $u_S$ and vector function $u_I$, respectively.  The equation for the latter can be represented as
\begin{equation}
\label{uI}
\partial_t u_I = \nabla \cdot ( D \nabla u_I ) + F(x)u_I - V(x) u_I
\end{equation}
where
$$u_I(t,x) = \left [\begin{array}{c} I(t,x) \\ V(t,x) \end{array} \right ], \quad D = \mathrm{diag}(D_I, D_V),$$
and
$$F(x) = \left[ \begin{array}{cc}
0 & kT_\infty(x) \\
0 & 0  \end{array} \right], \quad
V(x) = \left[ \begin{array}{cc}
\mu_I & 0 \\
-N\mu_I & \mu_V  \end{array} \right].$$
Using the properties of $T_\infty$ guaranteed by Theorem \ref{Tinfprop}, we note that the system in this form satisfies assumptions (A$1$)-(A$6$) of \cite{WangZhao}. 
Then, \cite[Theorem 3.1]{WangZhao} directly implies the conclusion under the assumption $\eta_0 < 0$, since this eigenvalue is exactly the spectral bound for the operator on the right side of \eqref{uI}.  Therefore, $E_c$ is a locally asymptotically stable equilibrium for \eqref{eq:system}.

To prove the second conclusion, we will show that the eigenvalue problem
\begin{equation}
\label{eigmatrix}
\nabla \cdot ( D \nabla \phi) + M(x) \phi = \eta \phi
\end{equation}
with $D$ and $M$ defined by \eqref{DM} possesses a positive eigenvalue with corresponding positive eigenvector under the assumption $\eta_0 > 0$.
To this end, we will again utilize a result from \cite{WangZhao}.  
From \eqref{DM}, we see that $M(x)$ is cooperative for all $x \in \Omega$ and $M(x_0)$ is irreducible for some $x_0 \in \Omega$ since, by Theorem \ref{Tinfprop}, $T_\infty(x) \geq 0$ for every $x \in \Omega$ and $T_\infty \not\equiv 0$.
Hence, by \cite[Theorem 2.2]{WangZhao}, $\eta_0$ as defined by \eqref{eta0} is an algebraically simple eigenvalue of \eqref{eigmatrix} with a strongly positive eigenvector satisfying
$$\mathrm{Re}(\eta) < \eta_0 \qquad \ \mathrm{for \ all} \ \eta \in \sigma \biggl (\nabla \cdot (D\nabla) + M \biggr )\setminus \{\eta_0\}.$$
Therefore, $\eta_0$ is a positive eigenvalue with positive eigenfunction, and the instability of $E_c$ follows immediately.
Finally, we note that that these results are valid for either boundary condition \eqref{BC1} or \eqref{BC2} as mentioned in \cite[Remarks 2.2 and 3.1]{WangZhao}.
\end{proof}

With the local stability proof complete, we turn to the proofs of the global asymptotic behavior of the system.

\begin{proof}[Proof (Theorem \ref{Tinfty})]
Similar to the proof of Theorem \ref{unique}, we define $u(t,x) = T(t,x) - T_\infty(x)$ so that $u$ satisfies
$$\partial_t u - D_T \Delta u = -\mu_T u - kTV.$$
Using the integrating factor $e^{\mu_T t}$ and recalling the positivity of $T(t,x)$ and $V(t,x)$ guaranteed by Theorem \ref{PosBounds}, this becomes
$$(\partial_t - D_T \Delta)[e^{\mu_T t} u] = - ke^{\mu_T t}TV \leq 0.$$
Invoking Corollary \ref{heatbound2} with $g \equiv 0$ yields
$$\Vert e^{\mu_T t} u(t) \Vert_\infty \leq \Vert u(0) \Vert_\infty$$
which is equivalent to
$$\Vert T(t) - T_\infty \Vert_\infty \leq \Vert T_0 - T_\infty \Vert_\infty e^{-\mu_T t}$$
for every $t \geq 0$, and the proof is complete.
\end{proof}

\begin{proof}[Proof (Theorem \ref{Asymp})]
We begin by writing the system in non-dimensionalized form.  Although one can arrive at the same results without doing so, using the dimensionless system will simplify the proof by reducing the number of parameters which appear and the complexity of related expressions.  Only within this proof will the non-dimensionalized version of the system be utilized.

First, define the dimensionless concentrations by
$$T^*(t,x) = \frac{T(t,x)}{T_c}, \quad I^*(t,x) = \frac{I(t,x)}{I_c}, \quad V^*(t,x) = \frac{V(t,x)}{V_c}$$
where $T_c, I_c$, and $V_c$ are constants to be determined.
Additionally, we scale the spatial and time dimensions by letting $$t^* = \frac{t}{t_c}, \qquad x^*= \frac{x}{x_c}.$$
Here the scaling constants will be determined so as to minimize the dimension of the resulting parameter space.  Of course, differing spatial parameters (i.e., one in each spatial component) could be introduced, but we assume a uniform directional scaling. Substituting these expressions within \eqref{eq:system}, we find
\begin{equation} 
\label{rescale}
\left.
\begin{aligned}
\frac{\partial T^*}{\partial t^*} - \frac{D_T t_c}{x_c^2} \Delta T^* \ &= \frac{t_c}{T_c}\lambda(x)  - \mu_Tt_c T^* - kt_cV_c T^*V^* \\
\frac{\partial I^*}{\partial t^*} - \frac{D_I t_c}{x_c^2}\Delta I^* &= k\frac{t_cT_cV_c}{I_c} T^*V^* - \mu_It_cI^* \\ 
\frac{\partial V^*}{\partial t^*} - \frac{D_V t_c}{x_c^2}\Delta V^* &= N \mu_I \frac{t_cI_c}{V_c} I^*  - \mu_Vt_c V^*. 
\end{aligned}
\right \}
\end{equation}
We fix the time and spatial scales using the $T$-cell diffusion and decay rates by choosing $t_c = \frac{1}{\mu_T}$ and $x_c = \sqrt{\frac{D_T}{\mu_T}}$.  Next, we choose the scaling for the dependent variables so as to eliminate parameters in each equation. In particular, this is effective when
$$T_c = \frac{\mu_V}{kN}, \qquad I_c= \frac{\mu_V \mu_T}{kN \mu_I}, \qquad V_c = \frac{\mu_T}{k}.$$
Dropping the starred notation, this finally yields the non-dimensionalized system
\begin{equation} 
\label{3CMPDE}
\left.
\begin{aligned}
\frac{\partial T}{\partial t} - \Delta T  &=  q(x)  - T - TV \\
\frac{\partial I}{\partial t} - \beta_1 \Delta I  &= \alpha_1 (TV - I) \\ 
\frac{\partial V}{\partial t} - \beta_2 \Delta V  &= \alpha_2 (I  - V)
\end{aligned}
\right \}
\end{equation}
where
$$\alpha_1 = \frac{\mu_I}{\mu_T}, \qquad \alpha_2 = \frac{\mu_V}{\mu_T}, \qquad \beta_1 = \frac{D_I}{D_T}, \qquad \beta_2 = \frac{D_V}{D_T}$$
and 
$$q(x) := \frac{kN}{\mu_T \mu_V} \lambda(x).$$
We note that the system now possesses only five parameters, and as we will see, the $\alpha$ terms control its rates of convergence to equilibrium, while the $\beta$ coefficients describe the rates of relative diffusion.
Finally, we rescale the steady state $T_\infty(x)$, defined by \eqref{Tinf}, using the same values of $x_c$ and $T_c$ so that, again dropping the starred notation, this function satisfies the linear, stationary PDE
\begin{equation}
\label{Tinf2} 
- \Delta T_\infty =  q(x)  - T_\infty.
\end{equation}
Due to the introduction of the scaling parameter $T_c$, our initial assumption of $\Vert R_0 \Vert_\infty < 1$ then implies the newly scaled equilibrium T-cell count satisfies $\Vert T_\infty \Vert_\infty < 1$.
Recall that all original parameters are assumed to be positive, and hence all new parameters retain this property.
To prove the result, then, we will show that the conclusion holds for solutions of \eqref{3CMPDE} as the exponential decay of the quantities $I$ and $V$ in the scaled space-time variables will imply the result of the theorem in the unscaled variables with a change in the constants $C_0$ and $a$.  

With the dimensionless equations in place, estimates for $T$ are immediate.  From Theorem \ref{Tinfty} we rescale the resulting inequality (or equivalently, perform the same steps within the proof on \eqref{3CMPDE} and \eqref{Tinf2}) to find
$$\Vert T (t) \Vert_\infty \leq \Vert T_\infty \Vert_\infty + \Vert T(0) - T_\infty \Vert_\infty  e^{-t} =: P(t).$$
As $\Vert T_\infty \Vert_\infty < 1$, let us define the quantities $M = \frac{1}{2}(1 + \Vert T_\infty \Vert_\infty)$ and $\tilde{M} =  \frac{1}{2}(1 + M)$ so that $\Vert T_\infty \Vert_\infty < M < \tilde{M} < 1$.  
Since $P(t)$ is decreasing, we further choose $$\tau_0 = \max \left \{ 0, \ln \left (\frac{\Vert T_0 - T_\infty\Vert_\infty}{M - \Vert T_\infty \Vert_\infty} \right ) \right \}$$
so that $P(t) \leq M$ for all $t \geq \tau_0$.

Next, let $a^* = \frac{1}{2}\min\{1, \alpha_1, \alpha_2\}$ and choose $a \in (0,a^*)$ small enough such that 
$$b := \tilde{M} - \frac{\alpha_1\alpha_2}{(\alpha_1 - a)(\alpha_2 - a)} M$$
satisfies $b > 0$. This choice of $a$ will allow us to take the constant in the exponential decay sufficiently large so that $T$ is uniformly bounded and the decay is preserved in time.
In particular, let $C_1 > 0$ satisfy $$C_1 > \max \left \{\frac{\sup_{s \in [0,\tau_0]} \Vert V(s) \Vert_\infty}{\Vert V(0) \Vert_\infty}, \frac{\sup_{s \in [0,\tau_0]} \Vert I(s) \Vert_\infty}{\Vert I(0) \Vert_\infty}\right \} e^{a\tau_0}$$
and
\begin{equation}
\label{cond}
C_1 \geq \frac{1}{b}\left ( 1 + \frac{\alpha_2\Vert I(0) \Vert_\infty}{(\alpha_2-a)\Vert V(0) \Vert_\infty} \right ).
\end{equation}

With these constants in place, notice that for every $t \in [0,\tau_0]$, we have
\begin{equation}
\label{decay1}
\Vert V(t) \Vert_\infty \leq \sup_{s \in [0,\tau_0]} \Vert V(s) \Vert_\infty
 < C_1 \Vert V(0) \Vert_\infty e^{-a\tau_0} 
\leq C_1 \Vert V(0) \Vert_\infty e^{-at}.
\end{equation}
Thus, the definition of $C_1$ implies that $V$ satisfies the decay estimate on the bounded interval $[0,\tau_0]$, and the same inequality holds for $I$.  It remains to prove these estimates for $t \geq \tau_0$.

Now, because solutions are sufficiently regular we may invoke Theorem \ref{LOR} and the Sobolev Embedding Theorem to deduce that the mapping $t \to \Vert V(t) \Vert_\infty$ is continuous.  Hence, we will utilize a continuity argument involving this norm in order to complete the proof.
Let $$\tau_1 = \sup\{ t > 0 : \Vert V(s) \Vert_\infty \leq C_1 \Vert V(0) \Vert_\infty e^{-as} \ \mbox{for \ every} \ s \in [0,t]\}$$
and notice that $\tau_1 > \tau_0$.
Then, for $t \in [\tau_0, \tau_1]$, we have both 
\begin{equation}
\label{VTest}
\Vert V(t) \Vert_\infty \leq C_1 \Vert V(0) \Vert_\infty e^{-at}
\qquad \mbox{and} \qquad  
\Vert T(t) \Vert_\infty \leq M.
\end{equation}
Deriving estimates on $I(t,x)$ as for its dimensionalized version, we see
$$\left (\frac{\partial}{\partial t} - \beta_1 \Delta \right ) \left [ e^{\alpha_1 t} I \right ]  = \alpha_1e^{\alpha_1 t}TV.$$
Integrating and using both Corollary \ref{heatbound1} and \eqref{VTest} implies
\begin{eqnarray*}
\Vert I(t) \Vert_\infty & \leq & \Vert I(0) \Vert_\infty e^{-\alpha_1 t} + \alpha_1 M e^{-\alpha_1 t} \int_0^t C_1 \Vert V(0) \Vert_\infty e^{(\alpha_1-a)s} \ ds\\
& \leq &  \Vert I(0) \Vert_\infty e^{-\alpha_1t} + \frac{\alpha_1}{\alpha_1 - a} \Vert V(0) \Vert_\infty  C_1 M \left ( e^{-at} - e^{-\alpha_1t} \right )\\
& \leq & C_2 e^{-at}
\end{eqnarray*}
where $$C_2 =  \Vert I(0) \Vert_\infty + \frac{\alpha_1}{\alpha_1 - a} \Vert V(0) \Vert_\infty  C_1 M.$$
Doing the same for $V$ and using the newly-derived estimate on $\Vert I(t) \Vert_\infty$, we find
$$\left (\frac{\partial}{\partial t} - \beta_2 \Delta \right ) \left [ e^{\alpha_2 t} V \right ]  = \alpha_2e^{\alpha_2 t}I$$
and thus
\begin{eqnarray*}
\Vert V(t) \Vert_\infty & \leq & \Vert V(0) \Vert_\infty e^{-\alpha_1 t} + \alpha_2 C_2 e^{-\alpha_2t} \int_0^t e^{(\alpha_2-a)s} \ ds\\
& \leq &  \Vert V(0) \Vert_\infty e^{-\alpha_2 t} + \frac{\alpha_2}{\alpha_2 - a}C_2 \left ( e^{-at} - e^{-\alpha_2 t} \right )\\
& \leq & \left (  \Vert V(0) \Vert_\infty  +   \frac{\alpha_2}{\alpha_2 - a} C_2 \right )  e^{-at}.
\end{eqnarray*}
Thus, the exponential decay of $\Vert V (t) \Vert_\infty$ continues on this time interval and we merely require the constant within this inequality to be strictly dominated by $C_1 \Vert V(0) \Vert_\infty$ to complete the argument.  After a brief calculation, we see that \eqref{cond} implies 
$$\Vert V(0) \Vert_\infty  + \frac{\alpha_2}{\alpha_2 - a} C_2 \leq \tilde{M} C_1\Vert V(0) \Vert_\infty,$$ and combining this with the estimate above, we find 
$$\Vert V(t) \Vert_\infty \leq \tilde{M} C_1 \Vert V(0) \Vert_\infty e^{-at}$$
for any $t \in [\tau_0, \tau_1]$.  Since $\tilde{M} < 1$, we see that $\tau_1$ cannot be finite, as this would contradict its definition as the supremum, and thus $\tau_1 = \infty$. Therefore, using \eqref{decay1} the estimate 
$$\Vert V(t)\Vert_\infty \leq C_1 \Vert V(0) \Vert_\infty e^{-at}$$ holds for all $t \geq 0$.
Since we have $\Vert I(t) \Vert_\infty \leq C_2 e^{-at}$ on the same interval, the exponential decay of $\Vert I(t) \Vert_\infty$ follows with the same rate.
Taking $C_0 = \max \{ C_1, C_2\}$ yields a uniform decay estimate on both non-dimensionalized quantities, namely
$$ \infnorm{I(t)} + \infnorm{V(t)}\leq C_0e^{-at}.$$
Finally, upon rescaling the non-dimensionalized system, the result holds for the original concentrations $I$ and $V$ with a change in the constants $C_0$ and $a$.

\end{proof}

We end this section with the proof of the associated corollary.

\begin{proof}[Proof (Corollary \ref{global})]
Assuming the condition on the parameters holds, we note that
$\Vert T_\infty \Vert_\infty \leq \frac{\Vert \lambda \Vert_\infty}{\mu_T}$ by Theorem \ref{Tinfprop}.  Hence, $\frac{Nk\Vert \lambda \Vert_\infty}{\mu_T \mu_V} < 1$ implies 
$$\Vert R_0 \Vert_\infty = \frac{Nk \Vert T_\infty \Vert_\infty}{\mu_V}< 1,$$ and the assumptions of Theorem \ref{Asymp} are satisfied, which implies the result.
\end{proof}

\section{Proofs of Lemmas}

Finally, to complete the paper, we include the proofs of lemmas from the previous sections.  The first few results are fairly straightforward applications of the ideas inherent within the Maximum Principle, but we include them for completeness.

\begin{proof}[Proof (Lemma \ref{Pos})]
While both results are classical, we will prove the former assertion and direct the reader to \cite{Evans} for details regarding the latter.
Define the positive and negative parts of $u$ by 
$$\begin{gathered}
u_+(t,x) = \max\{0, u(t,x) \} \\
u_-(t,x) = -\min\{0,u(t,x)\}
\end{gathered}$$ 
and notice that $u = u_+ - u_-$. With this in mind, we multiply \eqref{eq:scal_heat} by $u_-$ and integrate in both time and space. Then the left side is given by 
$$L := \int^s_0 \int_\W u_- \partial_t u \, dx dt - D \int^s_0 \int_\W u_- \D u \, dx dt \defeq  I + II, $$ 
where $s > 0$. Now, define $\W_s^- = \{ (\tau,x) \in (0,s] \times \W : u(\tau,x) \le 0 \}$, then $u_- = 0$ outside of $\W_s^-$ so we see $$I =  \int_{\W_s^-} u_- \partial_t u \, dxdt. $$ However, on this set $u = -u_-$ so $$I = -\int_{\W_s^-} u_- \partial_t u_- \, dxdt = - \frac{1}{2} \int_{\W_s^-} \partial_t \left(u_-^2 \right) \, dxdt.$$ However, since $u_- = 0$ outside of $\W_s^-$, integrating over $\W_s^-$ is the same as integrating over $(0,s]  \times  \W$. Thus $$I = - \frac{1}{2} \int_\W \int^s_0 \partial_t \left(u_-^2 \right) \, dt dx = -\frac{1}{2} \int_\W \left[u_-(s,x)^2 - u_-(0,x)^2 \right] dx. $$ Finally, $u_0 \ge 0$ implies that $u_-(0,x) \equiv 0$ so we find $I = - \tfrac{1}{2} \twonorm{u_-(s)}^2 \le 0.$

We use a similar set of steps for $II$: \begin{align*}
II &= -  D \int^s_0 \int_\W u_- \D u \, dxdt \\ &= D \int_{\W^-_s} u_- \D u_- \, dxdt \\ &= - D \int_{\W^-_s} \nabla u_- \cdot \nabla u_- \, dxdt + D \left. \int_0^s u_- \partder{u_-}{n} \right|_{\partial \W}\, dt.
\end{align*} Enforcing the homogeneous boundary condition, either \eqref{BC1} or \eqref{BC2}, yields $$II = - D \int^s_0 \int_\W \left|\nabla u_-(t,x) \right|^2 \, dxdt \le 0. $$ Thus we have $$L = I + II \le 0. $$ 

Next, considering the right side under the same operations, we have $$R := \int^s_0 \int_{\W}g(t,x) u_-(t,x) \, dxdt.$$ But both $g$ and $u_-$ are nonnegative so $R \ge 0$. 

Because $L=R$ with $L \leq 0$ and $R \geq 0$, we find that both must be zero.  Additionally, $L = 0$ forces $I = II = 0$, and we conclude $$\twonorm{u_-(s)} = 0.$$ This is only possible if $u_-(s,x) \equiv 0$. However, $s$ was an arbitrary element of $(0,\infty)$ and so $u_-(t,x) \equiv 0$ for all $x \in \W, t \in (0,\infty)$. But if $u_- = 0$, then $u = u_+ \ge 0$ which completes the proof.
\end{proof}

\begin{proof}[Proof (Corollaries \ref{heatbound1} and \ref{heatbound2})]
Define the function $$v(t,x) = \left( \infnorm{u_0} + \int^t_0 \infnorm{g(\tau)} d\tau \right)- u(t,x).$$ Then we notice that $$(\heat ) v = \infnorm{g(t)} - (\heat ) u = \infnorm{g(t)} - g(t,x) \geq 0.$$ 
Also $$v(0,x) = \infnorm{u_0} - u(0,x) = \infnorm{u_0} - u_0(x) \geq0.$$ Thus, by Lemma \ref{Pos}, $v(t,x) \ge 0$ for $x\in \W, t \geq 0$, and 
$$u(t,x) \le  \infnorm{u_0} + \int^t_0 \infnorm{g(\tau)} d\tau,$$ 
for all $x\in \W, t \geq 0.$ 
Taking the supremum over $x\in\W$ yields 
$$\infnorm{u(t)} \le  \infnorm{u_0} + \int^t_0 \infnorm{g(\tau)} d\tau,$$ for all $t \geq 0$ which completes the proof.
The same method then applies to prove Corollary \ref{heatbound2}.
\end{proof}



\medskip
Received xxxx 20xx; revised xxxx 20xx.
\medskip

\end{document}